\documentclass[12pt]{amsart}

\usepackage{amsthm,amscd,amssymb,amsmath,tikz-cd,mathtools,enumerate}
\usepackage[mathscr]{eucal}
\usepackage{textcomp}

\theoremstyle{defn}
\newtheorem{defn}{Definition}[section]

\theoremstyle{plain}

\newtheorem{theorem}[defn]{Theorem}
\newtheorem{proposition}[defn]{Proposition}

\newtheorem{conjecture}[defn]{Conjecture}

\theoremstyle{remark}
\newtheorem{remark}[defn]{Remark}
\newtheorem{example}[defn]{Example}

\newtheorem{question}[defn]{Question}

\newcommand{\Gr}{\mathrm{Gr}}

\usepackage[margin=0.9in]{geometry}
\begin{document}

\title{{\bf P=W Phenomena}}

\author{Andrew Harder, Ludmil Katzarkov, Victor Przyjalkowski}

\address{\emph{Andrew Harder}
\newline
\textnormal{Department of Mathematics, Lehigh University, Christmas-Saucon Hall, 14 E. Packer Ave., Bethlehem, PA, USA, 18015.}
\newline
\textnormal{\texttt{anh318@lehigh.edu}}}

\address{\emph{Ludmil Katzarkov}
\newline
\textnormal{University of Miami, Coral Gables, FL, 33146, US}
\newline
\textnormal{HSE University, Russian Federation, Laboratory of Mirror Symmetry, NRU HSE, 6 Usacheva str., Moscow, Russia, 119048.}
\newline
\textnormal{\texttt{lkatzarkov@gmail.com}}}

\address{\emph{Victor Przyjalkowski}
\newline
\textnormal{Steklov Mathematical Institute of Russian Academy of Sciences, 8 Gubkina street, Moscow 119991, Russia.}
\newline
\textnormal{HSE University, Russian Federation, Laboratory of Mirror Symmetry, NRU HSE, 6 Usacheva str., Moscow, Russia, 119048.}
\newline
\textnormal{\texttt{victorprz@mi.ras.ru, victorprz@gmail.com}}}

\thanks{A.\,Harder was supported during part of this work by the Simons Collaboration in Homological Mirror Symmetry. L.\,Katzarkov was supported by Simons research grant, NSF DMS 150908, ERC Gemis, DMS-1265230, DMS-1201475, OISE-1242272 PASI, Simons collaborative Grant --- HMS,
Simons investigator grant --- HMS; he is partially supported by Laboratory of Mirror Symmetry NRU HSE, RF Government grant, ag. \textnumero~14.641.31.0001.
V.\,Przyjalkowski was partially supported by Laboratory of Mirror Symmetry NRU HSE, RF Government grant, ag. \textnumero~14.641.31.0001.
He is Young Russian Mathematics award winner and would like to thank its sponsors and jury.}

\maketitle

\begin{abstract}
In this paper, we describe recent work towards the mirror P=W conjecture, which relates the weight filtration on a cohomology of a log Calabi--Yau manifold to the perverse Leray filtration on the cohomology of the homological mirror dual log Calabi--Yau manifold, taken with respect to the affinization map. This conjecture extends the classical relationship between Hodge numbers of mirror dual compact Calabi--Yau manifolds, incorporating tools and ideas which appear in the fascinating and groundbreaking works of de Cataldo, Hausel, and Migliorini \cite{dcmh}, and de Cataldo and Migliorini \cite{dcm}. We give a broad overview of the motivation for this conjecture, recent results towards it, and describe how this result might arise from the SYZ formulation of mirror symmetry. This interpretation of the mirror P=W conjecture provides a possible bridge between the mirror P=W conjecture and the well-known P=W conjecture in nonabelian Hodge theory.
\end{abstract}

\section{Introduction}

The P=W conjecture was introduced in the seminal work
\cite{dcmh}. In this paper we introduce  a new read of
the P=W conjecture in the case of Landau--Ginzburg models.

Throughout this article, the reader may interpret our use of the phrase ``mirror symmetry''
to refer to homological mirror symmetry unless otherwise specified. We will abuse convention
and neglect specific complex and symplectic structures because our goal is to discuss relations
between Hodge-theoretic invariants which are largely insensitive to such choices.

Given a smooth quasiprojective variety $U$, an interesting combinatorial invariant is its dual intersection complex. Choosing a projective simple normal crossings compactification, $X$ of $U$  with $D = X \setminus U$, we let $\Gamma(D)$ be the dual intersection complex of $D$. The homotopy type of $\Gamma(D)$ is an invariant of $U$ (see \cite[Corollary 3]{KS06}), and the homology of $\Gamma(D)$ determines the weight $0$ part of the Deligne's canonical mixed Hodge structure on compactly supported cohomology
of $U$.

Often, in mirror symmetry, one is interested in pairs $(X,D)$ where $D$ is a simple normal crossings anticanonical divisor in $X$. We call such a pair {\it log Calabi--Yau}. We abuse notation and call $U = X \setminus D$ log Calabi--Yau if the pair $(X,D)$ is.

In \cite{ghk}, Gross, Hacking, and Keel start with a pair $(X,D)$ where $X$ is a surface and~$D$ is an anticanonical cycle of rational curves and construct a likely two-dimensional mir\-ror~$U^\vee$ to~$\mbox{$U = X \setminus D$}$, as the spectrum of a certain ring of functions. Therefore, the ring of functions of~$U^\vee$, tautologically, has dimension 2.
On the other hand, according to Auroux~\cite{ADC},  if~$(X,E)$ is a pair consisting of a del Pezzo surface $X$ and a smooth anticanonical divisor,  then the SYZ mirror to $U = X \setminus E$ is a rational elliptic surface with a fiber removed. Hence the spectrum of its ring of functions is of dimension 1. So we observe that if $(X,D)$ is a pair made up of a rational surface $X$ and a reduced simple normal crossings anticanonical divisor~$D$, then the dimension of the cone over the dual intersection complex of $U = X \setminus D$ equal to the dimension of~$\mathrm{Spec}(H^0(U^\vee, \mathscr{O}_{U^\vee}))$, where $U^\vee$ is the mirror of $U$.

Generally, one might expect that some analogue of this holds in arbitrary dimension, but as
it stands, this is not yet a precise question.
The goal of the work described in this article is to expand on this observation and to make precise conjectures in terms of the cohomology rings and mixed Hodge structures on mirror pairs of log Calabi--Yau manifolds $U$ and $U^\vee$.

In the following, all cohomology groups are taken with complex coefficients for the sake of simplicity. If $U$ and $U^\vee$ are mirror log Calabi--Yau manifolds, then we expect at first approximation that $H^*_c(U)$ and $H^*_c(U^\vee)$ are isomorphic as vector spaces (see \cite[Table 1]{KKP17}) with different gradings. By Poincar\'e duality, $H^i_c(U) \cong H^{\dim U - i}(U)$, hence we may equivalently deal with the cohomology rings of $U$ and $U^\vee$. Both $H^*(U)$ and $H^*(U^\vee)$ admit a mixed Hodge structure, which is composed of a decreasing Hodge filtration $F^\bullet$ and an increasing weight filtration~$W_\bullet$. We will define
\[
h^{p,q}(U) = \dim \Gr_F^q H^{p+q}(U).
\]
In analogy with classical mirror symmetry for compact Calabi--Yau varieties, we might expect that if $U$ and $U^\vee$ are a homological mirror pair of log Calabi--Yau manifolds of dimension $d$, then
\begin{equation}\label{eq:logcyhn}
h^{p,q}(U) = h^{d - p,q}(U^\vee).
\end{equation}
This seems to be true --- it is checked in many cases in \cite{HKP} --- but it ignores the weight filtration in cohomology. It would be desirable to determine whether the weight filtration on~$H^*(U)$ is reflected by a filtration on the cohomology of $U^\vee$.
The first step in this is to remark that the geometry of $D = X \setminus U$ and the residues of holomorphic forms on $X$ with log poles along $D$  can be used to determine the weight filtration $W_\bullet$ on $H^*(U)$. The weight filtration depends on the existence of a projective simple normal crossings compactification $X$ of $U$, but is independent of the choice of compactification, hence it is a canonical invariant of $U$. So if a mirror dual filtration exists, is plausible that it can be constructed via information dual to that of the components of $D$ but be independent the choice of $D$.

Starting with a log Calabi--Yau manifold $U$ and a simple normal crossings com\-pac\-ti\-fi\-ca\-ti\-on~$X$ of $U$ with $D = X \setminus U$, each irreducible component $D_i$, $i=1,\ldots,k$, of $D$ determines a regular function $w_i$ on the mirror $U^\vee$, see~\cite{aur1, aur2, aak}. Therefore, if there is a filtration on $H^*(U^\vee)$ dual to the weight filtration on $H^*(U)$, it should be determined by the func\-ti\-ons~$w_{1},\ldots , w_{k}$.

There are several possible filtrations on cohomology that can be constructed from $(w_1,\dots, w_k)$, but the most relevant seems to be the {\it flag filtration} \cite{dcm}, which is defined as follows. Let $w$ denote the map $(w_1,\dots, w_k)\colon U^\vee\rightarrow \mathbb{C}^k$. Choose a generic flag of  linear sub\-spa\-ces~$\Lambda_k \subset \Lambda_{k-1} \subset \dots  \subset \Lambda_0 =\mathbb{C}^k$ so that $\dim \Lambda_i = k-i$ and let $U^\vee_i = w^{-1}(\Lambda_i)$. Then, for any coefficient ring $R$, the flag filtration on $H^*(U^\vee;R)$ is defined as\footnote{Note that this agrees with the definition of \cite{dcm} up to a shift by $j$.}
\[
P_{r}H^j(U^\vee;R) = \ker(H^j(U^\vee;R) \longrightarrow H^j(U^\vee_{r+1};R)).
\]
According to de Cataldo and Migliorini \cite{dcm}, if $w$ is proper, then $P_\bullet$
can be identified with the
\emph{perverse Leray filtration} of the map $w$, hence it only depends on the map $w$. In all of the cases
that we know of, the maps $w_{1},\dots, w_{k}$ generate $\mathbb{C}[U^\vee]$, in which case the map $w: U^\vee \rightarrow \mathrm{im}(U^\vee)$ is the affinization map of $U^\vee$, hence, in these cases at least, $P_\bullet$ is {\it intrinsic to $U^\vee$} and does not depend on our original
choice of $w_{1},\dots, w_{k}$. Thus we have two filtrations, which are built from data which correspond to one another under mirror symmetry, and which are intrinsic to $U$ and $U^\vee$ respectively.
\begin{defn}
Consider a quasiprojective variety $M$ over the complex numbers $\mathbb{C}$ %, assume that its af\-fi\-ni\-za\-ti\-on~$\mbox{$\mathrm{Spec}(\mathbb{C}[M])$}$ is an
%affine variety, and, for simplicity, 
and assume that the affinization map $f^\mathrm{aff}\colon M \rightarrow \mathrm{Spec}(\mathbb{C}[M])$ is proper. We define the \emph{perverse mixed Hodge polynomial} of a quasiprojective variety $M$ to be
\[
{\mathrm{PW}}_{M}(u,t,w,p) = \sum_{a,b,r,s} (\dim \Gr_F^a \Gr^W_{s+b}\Gr^P_{r} (H^{s}(M)))u^a t^{s}w^{b}p^r,
\]
where $P_\bullet$ is the flag filtration taken with respect to $f^\mathrm{aff}$ and $W$ denotes the $\mathbb{C}$-linear extension of the weight filtration.
\end{defn}
We make the following conjecture.
\begin{conjecture}[Mirror P=W conjecture]\label{con:filt}
Let $U$ be a log Calabi--Yau variety and assume that its homological mirror $U^\vee$ is also a log Calabi--Yau variety whose dimension is the same as that
of $U$. Let $d = \dim U = \dim U^\vee$. Then
\begin{equation}\label{eq:toptest}
{\mathrm{PW}}_{U}(u^{-1}t^{-2}, t,p, w)u^dt^d = {\mathrm{PW}}_{U^\vee}(u,t,w,p).
\end{equation}
\end{conjecture}
\begin{remark}
The maximal depth of the perverse Leray filtration on $U$ is equal to the dimension of~$\mathrm{Spec}(H^0(U,\mathscr{O}_U))$, so it corresponds to the maximal dimension of a cell in the dual intersection complex for $U^\vee$. Therefore, Conjecture \ref{con:filt} is consistent with the observations that we have made in the case of surfaces.
\end{remark}
\begin{remark}
One can show (\cite{HKP}) that if $U$ is affine, then Conjecture \ref{con:filt} implies that $U^\vee$ has Hodge--Tate cohomology ring. If $U^\vee$ is log Calabi--Yau, then having a Hodge--Tate cohomology ring implies (but is not equivalent to) the fact that the dimension of the dual intersection complex for $U^\vee$ is the same as the complex dimension of $U^\vee$.
\end{remark}

We would like to emphasize that the mirror P=W conjecture is intimately related to several conjectures that have already appeared in the literature. First, if $X$ is a Fano manifold, its mirror is a Landau--Ginzburg model $(Y,w)$. In fact, as explained by the second author, Kontsevich, and Pantev \cite{KKP17}, if $D$ is a simple normal crossings anticanonical divisor in $X$, the variety~$U^\vee$ is the homological mirror of $U = X \setminus D$ as above, and $w_1,\dots, w_k$ are functions on $U^\vee$ associated to the components of $D$, then the homological mirror to $X$ is the Landau--Ginzburg pair $(U^\vee, w_1 + \dots + w_k)$. Furthermore, there are two major conjectures made which relate the Hodge theory of $X$ and new Hodge theoretic invariants of $(U,w_1+ \dots + w_k)$. We will explain in Section \ref{sect:kkp} how these conjectures follow from the mirror P=W conjecture in the case where $D$ is smooth. Analogous results should hold when $D$ is no longer smooth, however we believe that further study of the Hodge theory developed in \cite{KKP17} is necessary to prove this.

Additionally, as the name suggests, the mirror P=W conjecture is related (in spirit at least) to the P=W conjecture of de Cataldo, Hausel, and Migliorini \cite{dcmh}. We recall that if $\mathcal{M}_\mathrm{B}$ and~$\mathcal{M}_\mathrm{H}$ denote Betti and Higgs moduli spaces associated to certain geometric and representation theoretic data, then $\mathcal{M}_\mathrm{H}$ admits a Hitchin map, $h\colon \mathcal{M}_\mathrm{H} \rightarrow \mathbb{C}^{d}$, where $\dim_\mathbb{C} \mathcal{M}_\mathrm{H} = 2d$. It is known that $\mathcal{M}_\mathrm{B}$ and $\mathcal{M}_\mathrm{H}$ are diffeomorphic quasiprojective varieties which are not deformation equivalent. In \cite{dcmh}, de Cataldo, Hausel, and Migliorini conjecture (and prove in a large number of cases) that under a certain diffeomorphism,
\[
W_{2i} H^n(\mathcal{M}_\mathrm{B}) = W_{2i+1}H^n(\mathcal{M}_\mathrm{B}) = P_iH^n(\mathcal{M}_\mathrm{H})
\]
where the perverse Leray filtration on $H^n(\mathcal{M}_\mathrm{H})$ is taken with respect to the Hitchin map.

The Betti moduli space is expected to be log Calabi--Yau \cite{Simpson}, and admits a special Lagrangian torus fibration obtained from the Hitchin fibration. Therefore, a restatement of the classical P=W conjecture is that the weight filtration on the cohomology of $\mathcal{M}_\mathrm{B}$ can be obtained from the geometry of this given special Lagrangian torus fibration. We refer to this as the geometric P=W conjecture.

 In Section \ref{sect:hkkp}, we explain work in progress by the first two authors along with Pantev and Kontsevich how the mirror P=W conjecture might be obtained from the SYZ formulation of mirror symmetry. Therefore, we obtain a geometric version of the mirror P=W conjecture, in which the input data is reminiscent of the input data for the geometric P=W conjecture.

Once these geometric considerations are out of the way we will discuss possible relationships between the classical and mirror P=W conjectures.

\medskip

{\bf Acknowledgements.} The authors would like to thank Valery Lunts and Tony Pantev for enlightening conversations.

\section{Connections to the conjectures of \cite{KKP17}}\label{sect:kkp}

Let $X$ be a Fano manifold, then its mirror dual is expected to be a Landau--Ginzburg model $(Y,w)$, where~$Y$ is a log Calabi--Yau variety and $w\colon Y\rightarrow \mathbb{C}$ is a regular function. The authors of \cite{KKP17} constructed Hodge-theoretic invariants, $f^{p,q}(Y,w)$, of a Landau--Ginzburg model. We may define
\[
f^{p,q}(Y,w) = \dim \mathrm{Gr}^F_qH^{p+q}(Y,V)
\]
where $V$ is a smooth fiber of $w$ and $H^{p+q}(Y,V)$ is equipped with the natural mixed Hodge structure on relative cohomology. This is equivalent to the definition of \cite{KKP17} by \cite{Ha17}. The second author, along with Kontsevich and Pantev has conjectured that if $X$ and $(Y,w)$ form a homological mirror pair, then
\begin{equation}\label{eq:kkpeqq}
f^{p,q}(Y,w) = h^{\dim X - p,q}(X).
\end{equation}
When the map $w$ is proper we expect that $X$ admits a smooth anticanonical divisor $D$ so that~$X \setminus D$ and $Y$ form a homological mirror pair. Therefore, Equality~\eqref{eq:logcyhn} should hold between $Y$ and~$X \setminus D$. Furthermore, we expect that a general smooth fiber $V$ of $w$ is Calabi--Yau, and is the homological mirror of~$D$, so we expect that
\begin{equation}\label{eq:cptcyhn}
h^{p,q}(V) = h^{\dim X -1 -p,q}(D).
\end{equation}
We show that Conjecture \ref{con:filt} links Equality~\eqref{eq:kkpeqq} with Equalities~\eqref{eq:logcyhn} and \eqref{eq:cptcyhn}. Stated roughly, this takes the following form.
\begin{theorem}[\cite{HKP}]
Let $X$ be a projective manifold with a smooth anticanonical divisor $D$ in it, and let $U = X\setminus D$. Let $(Y,w)$ be a Landau--Ginzburg model so that $w$ is proper and let $V$ be a smooth fiber of $w$. If $V$ and $D$ satisfy (\ref{eq:cptcyhn}), and $Y$ and $U$ satisfy Equality~\eqref{eq:logcyhn}, then Conjecture \ref{con:filt} implies Equality~\ref{eq:kkpeqq}.
\end{theorem}
Finally, for a Landau--Ginzburg model $(Y,w)$ consider a monodromy action
\[
M \colon H^i(Y,V) \longrightarrow H^i(Y,V)
\]
obtained by letting $V$ vary in a small circle around $\infty$. Taking the logarithm of $M$, which we denote $N$, we obtain a filtration on $H^i(Y,V)$, which we will call $\mathrm{Mon}_\bullet$ (see \cite{LP18} for a precise definition). We let
\[
h^{p,q}(Y,w) = \dim \mathrm{Gr}^\mathrm{Mon}_{p} H^{q}(Y,V).
\]
The authors of \cite{KKP17} conjecture that if $X$ is a Fano manifold and $(Y,w)$ is its homological mirror, then
\begin{equation}\label{eq:fanoht}
h^{p,q}(Y,w) = f^{p,q}(Y,w).
\end{equation}
We have proved the following relationship between Equality \eqref{eq:fanoht} and Conjecture \ref{con:filt}.

\begin{theorem}[\cite{HKP}]
Let $X$ be a Fano manifold with a smooth anticanonical divisor $D$ in $X$, and $(Y,w)$ be a Landau--Ginzburg model so that $w$ is proper. Assume that Conjecture \ref{con:filt} holds between $Y$ and $U = X \setminus D$. Then
\[
f^{p,q}(Y,w) = h^{p,q}(Y,w).
\]
\end{theorem}

Now consider the case of dimensions $2$ and $3$. Landau--Ginzburg models (satisfying a half of Homological Mirror Symmetry) for del Pezzo surfaces
are constructed in \cite{AKO06}. According to loc. cit., the tame compactified mirror Landau--Ginzburg model for del Pezzo surface of degree~$d$ is a rational elliptic surface whose fiber
over infinity is a wheel of $d$ smooth rational curves. In particular, Equality~\eqref{eq:cptcyhn} holds for del Pezzo surfaces since for all of them the anticanonical linear system
is non-empty. Moreover, by~\cite{LP18}, Equalities~\eqref{eq:kkpeqq} and~\eqref{eq:fanoht} hold.
This implies that Conjecture~\ref{con:filt} holds for a pair $(X,D)$, where $X$ is a smooth del Pezzo surface and $D$ is a smooth anticanonical divisor on it.

Any smooth Fano threefold has \emph{toric Landau--Ginzburg model} (see~\cite{Prz18} for details), that is a Landau--Ginzburg model whose total space is an algebraic torus
$(\mathbb C^*)^3$ and which satisfies certain conditions. In~\cite{Prz17} their \emph{log Calabi--Yau compactifications} are constructed.
In particular, they are tame compactified Landau--Ginzburg models $(Z,f)$, whose fibers form an anticanonical linear system,
so fibers have trivial dualizing sheaves. The dual intersection complex of the fiber over infinity is homotopic to a sphere, and components of the fiber over infinity are smooth rational surfaces. Finally, one has~$h^{p,q}(Z)=0$ if $p\neq q$. That is, the cohomology ring of $Z$ is Hodge--Tate. Under these assumptions in~\cite{Ha17} it has been shown that the $f^{p,q}(Y,w)$-diamond for the compactified Landau--Ginzburg model is of the form
$$
\begin{matrix}
&&&0&&& \\
&&0&&0&&\\
&0&& k_Y &&0& \\
1\qquad&& ph-2 + h^{1,2}(Z) && ph-2 + h^{2,1}(Z) &&\qquad1 \\
&0&& k_Y &&0& \\
&&0&&0&&\\
&&&0&&&
\end{matrix}
$$
where
$$
ph=\dim\Bigg(\mathrm{coker}\Big(H^2\big(Z,\mathbb{R}\big)\to H^2\big(V,\mathbb{R}\big)\Big)\Bigg)
$$
and $k_Y$ is defined to be
$$
k_Y = \sum_{s\in \Sigma}(\rho_s-1),
$$
where $\Sigma$ is a set of critical values of $w$ and $\rho_s$ is the number of irreducible components of $w^{-1}(s)$,
see~\cite[Theorem 22]{Prz13} and~\cite[Conjecture 1.1]{PSh15a}.
Using this result, in~\cite{CP18} it is proven that Equality~\eqref{eq:kkpeqq} holds for Fano threefolds.
The second conjecture from~\cite{KKP17} for Fano threefolds (Equality~\eqref{eq:fanoht}) is given by the following result,
which is an extension of the main result of Shamoto in~\cite{Sh17}.
\begin{proposition}[\cite{HKP}]\label{corollary:lastones}
Let $(Y,w)$ be a Landau--Ginzburg model, and assume that $w$ is a proper map.
If $H^i(Y)$ is Hodge--Tate for all $i$, then
\[
h^{p,q}(Y,w) = f^{p,q}(Y,w)
\]
for all $p,q$.
\end{proposition}

The converse to the results presented above also holds for Fano varieties in dimension $2$ and $3$.

\begin{theorem}[\cite{HKP}]
Let $(X,D)$ be a pair consisting of smooth Fano surface or threefold $X$ and a smooth anticanonical divisor $D$ on it.
Let $(Y,w)$ be its compactified Landau--Ginzburg model constructed in~\cite{AKO06} and~\cite{Prz17}.
Then Conjecture~\ref{con:filt} holds for them.
\end{theorem}

Equalities~\eqref{eq:kkpeqq} and~\eqref{eq:fanoht} also hold for a smooth toric weak Fano threefolds $X_\Delta$ for which the map~$H^2(X_\Delta) \rightarrow H^2(D)$ is injective for $D$ a smooth anticanonical divisor, see~\cite{Ha17}. Equality~\eqref{eq:cptcyhn} holds in this case because both $D$ and the fibers of the Landau--Ginzburg mirror of $X_\Delta$ are K3 surfaces.

\section{Example}
Let $X$ be an intersection of two quadrics in $\mathbb P^5$. According to~\cite{Prz13}, its toric Landau--Ginzburg model is
$$
f=\frac{(x+1)^2(y+1)^2}{xyz}+z.
$$
In particular, it satisfies mirror symmetry conjecture of variations of Hodge structures (that relates Gromov--Witten theory of $X$
and periods of the family of fibers for $f$), it is related to a toric degeneration of $X$, and it admits
a compactification to a tame compactified Landau--Ginzburg model, see below.
We suggest the compactification as a homological mirror dual to $X$.

The toric Landau--Ginzburg model is obtained in the following way. We start with Givental's Landau--Ginzburg model for $X$, defined as the subset of $(\mathbb{C}^*)^5=\mathrm{Spec}(\mathbb{C}[x^{\pm 1}, y^{\pm 1}, z^{\pm 1}, u^{\pm 1}, v^{\pm 1}])$
\begin{equation}\label{eq:lgeq}
x + y = 1, \quad z + u = 1
\end{equation}
equipped with the superpotential
\[
w = v + \dfrac{1}{xyzuv}.
\]
The anticanonical divisor of $X$ is the intersection of $X$ with another quadric in $\mathbb{P}^5$. The monomials in the expression of $w$ given here correspond to hyperplane sections of $X$, hence the pair of maps
\[
w_1 = u, \qquad w_2 = \dfrac{1}{xyzuv}
\]
correspond to hyperplane sections of $X$. There is a birational map from $(\mathbb{C}^*)^3$ to the locus cut out by \eqref{eq:lgeq},
\[
(x', z',v) \mapsto \left( \dfrac{x'}{x' + 1}, \dfrac{1}{x'+1}, \dfrac{z'}{z'+1}, \dfrac{1}{z'+1}, v\right)
\]
under which $w_1$ and $w_2$ become
\[
w_1 = v, \quad w_2 = \dfrac{(x'+1)^2(z'+1)^2}{x'z'v}.
\]
In \cite{HKP}, we construct a fiberwise compactification of this map which we will call $Y^\circ$. We refer to the extensions of $f_1$ and $f_2$ to $Y^\circ$ as $w_1$ and $w_2$. In \cite{Prz13}, a relative compactification with respect to the map $f = f_1 + f_2$ is constructed, which we call $Y$. We refer to the extension of $f$ to $Y$ as $w$. There is a divisor $D_h$ in $Y$, which is biregular to $I_{10} \times \mathbb{C}$ (where $I_n$ denotes a cycle of $n$ rational curves), so that the map $w$ restricted to $D_h$ is projection onto the second factor, and $Y^\circ = Y \setminus D_h$.

In \cite{HKP}, following results of \cite{CP18}, we compute the mixed Hodge structure on the compactly supported cohomology of $Y$ to be:
\begin{itemize}
\item[-] $H^6_c(Y;\mathbb{Q}) \cong \mathbb{Q}(-3)$;
\item[-] $H^4_c(Y;\mathbb{Q}) \cong \mathbb{Q}(-2)^{21}$;
\item[-] $H^3_c(Y;\mathbb{Q}) \cong \mathbb{Q}(0)$;
\item[-] $H^2_c(Y;\mathbb{Q}) \cong \mathbb{Q}(-1)^3$.
\end{itemize}
All other cohomology groups vanish. This allows us to compute $H^*(Y^\circ;\mathbb{Q})$. We can compute that $H^*_c(D_h;\mathbb{Q})$ is given by:
\begin{itemize}
\item[-] $H^2_c(D_h;\mathbb{Q}) \cong \mathbb{Q}(-1)$;
\item[-] $H^3_c(D_h;\mathbb{Q}) \cong \mathbb{Q}(-1)$;
\item[-] $H^4_c(D_h;\mathbb{Q}) \cong \mathbb{Q}(-2)^{10}$.
\end{itemize}
Again, all other cohomology groups vanish. We can use this, along with the standard long exact sequence in mixed Hodge structures
\[
\dots \rightarrow H^i_c(Y^\circ;\mathbb{Q}) \rightarrow H^i_c(Y;\mathbb{Q}) \rightarrow H^i_c(D_h;\mathbb{Q}) \rightarrow \dots
\]
to compute the mixed Hodge structure on $H^i_c(Y^\circ)$, which is dual to that of $H^{6-i}(Y^\circ)(-3)$. We find that:
\begin{itemize}
\item[-] $H^0(Y^\circ;\mathbb{Q}) \cong \mathbb{Q}(0)$;
\item[-] $H^2(Y^\circ;\mathbb{Q})$ is an extension of $\mathbb{Q}(-2)$ by $\mathbb{Q}(-1)^{11}$;
\item[-] $H^3(Y^\circ;\mathbb{Q})\cong \mathbb{Q}(-3)$;
\item[-] $H^4(Y^\circ;\mathbb{Q}) \cong \mathbb{Q}(-2)^2$.
\end{itemize}
Choose a flag $H_2\subset H_1\subset H_0=\mathbb A^2$, where $H_2$ is a point, and $H_1$ is a general line inside of $\mathbb{C}^2$. The fiber $Y_2$ over a point compactifies to an elliptic curve. The preimage $Y_1$ of $H_1$ inside of $Y^\circ$ compactifies to a fiber of $w$ with an $I_{10}$ configuration of rational curves removed. This then allows us to compute the perverse Leray filtration on $H^*(Y)$ and its mixed Hodge structure.
\begin{itemize}
\item[-] The kernel of $H^2(Y^\circ;\mathbb{Q}) \rightarrow H^2(Y_1;\mathbb{Q})$ is isomorphic to $\mathbb{Q}(-2)^8$.
\item[-] The map $H^0(Y^\circ;\mathbb{Q}) \rightarrow H^0(Y_1;\mathbb{Q})$ is an isomorphism.
\item[-] For $i \neq 0, 2$, the map $H^i(Y^\circ;\mathbb{Q}) \rightarrow H^i(Y_1;\mathbb{Q})$ is trivial.
\end{itemize}
Therefore,
\[
P_4H^4(Y^\circ;\mathbb{Q}) \cong \mathbb{Q}(-2)^8, \quad P_3H^3(Y^\circ;\mathbb{Q}) \cong \mathbb{Q}(-3), \quad P_2H^2(Y^\circ;\mathbb{Q}) \cong \mathbb{Q}(-1)^2,
\]
and $P_iH^i(Y^\circ;\mathbb{Q}) = 0$ otherwise. This should correspond to the graded pieces of the cohomology of the homological mirror which are of the weight equal to the rank of the cohomology group containing them. It's also easy to compute the restriction maps $H^i(Y^\circ) \rightarrow H^i(Y_2)$.
\begin{itemize}
\item[-] The map $H^0(Y^\circ;\mathbb{Q}) \rightarrow H^0(Y_2;\mathbb{Q})$ is an isomorphism.
\item[-] The map $H^2(Y^\circ;\mathbb{Q}) \rightarrow H^2(Y_2;\mathbb{Q})$ has is surjective.
\end{itemize}
Therefore, $P_3H^2(Y^\circ;\mathbb{Q})$ is an extension of $\mathbb{Q}(-2)$ by $\mathbb{Q}(-1)^2$, and $P_{i+1}H^i(Y^\circ;\mathbb{Q}) \cong 0$ is zero if $i \neq 2$. Furthermore, we deduce that $P_2H^0(Y^\circ;\mathbb{Q}) \cong \mathbb{Q}(0)$, $P_2H^2(Y^\circ;\mathbb{Q})\cong \mathbb{Q}(-1)$, and~$P_{i+2}H^i(Y^\circ;\mathbb{Q}) \cong 0$ if $i \neq 0$. Therefore,
\[
\mathrm{PW}_{Y^\circ}(u,v,w,p) = u^3t^3w^3 + 2ut^2 + 2u^2t^4 + p(8ut^2 + u^2t^2w^2) + p^2(1 + ut^2).
\]
Now we compute the perverse mixed Hodge polynomial of the prospective mirror of $Y^\circ$. Remind that $X$ is the intersection of two generic quadrics inside of $\mathbb{P}^5$, and we let $D_1, D_2$ be hyperplane sections of $X$ chosen generically, and we let $E = D_1 \cap D_2$. We note that $E$ is a smooth elliptic curve, $D_1$ and $D_2$ are del Pezzo surfaces of degree 4, hence $H^2(D_1;\mathbb{Q}) \cong H^2(D_2;\mathbb{Q}) \cong \mathbb{Q}(-1)^5$. Finally, it is known that $H^{2i}(X;\mathbb{Q}) \cong \mathbb{Q}(-i)$ if $0 \leq i \leq 3$ and that $H^3(X;\mathbb{Q}) \cong H^1(C;\mathbb{Q})(-1)$ for a curve of genus 2, hence $h^{1,2}(X) = h^{2,1}(X) = 2$. Using the standard spectral sequence to compute the weight-graded pieces of the mixed Hodge structure on $U = X \setminus (D_1 \cup D_2)$ (see e.g. \cite[Proposition 8.3.4]{voisinI}), we see that
\[
\mathrm{Gr}^5_W H^3(U;\mathbb{Q}) \cong H^1(E;\mathbb{Q})(-2), \quad \mathrm{Gr}_W^4 H^3(U;\mathbb{Q}) \cong \mathbb{Q}(-2)^8, \quad \mathrm{Gr}_W^3H^3(U;\mathbb{Q}) \cong H^3(X;\mathbb{Q}),
\]
\[
\mathrm{Gr}^2_WH^1(U;\mathbb{Q}) \cong \mathbb{Q}(-1), \quad \mathrm{Gr}^0_WH^0(U;\mathbb{Q}) \cong \mathbb{Q}(0),
\]
and all other weight-graded pieces in cohomology vanish. Finally, $U$ is an affine variety, so the perverse Leray filtration with respect to the identity map can be computed to be precisely:
\[
P_{3-i}H^i(U;\mathbb{Q}) \cong \begin{cases}  H^i(U;\mathbb{Q}) & \text{ if } j = i \\
 0 & \text{ otherwise.}
\end{cases}
\]
Therefore the perverse mixed Hodge polynomial of $U$ is
\[
\mathrm{PW}_U(u,v,w,p) = p^3 + p^2utw + t^3((u^3+ u^2)w^2 + 8u^2w + (2u^2 + 2u)).
\]
We can check:
\begin{align*}
\mathrm{PW}_U(u^{-1}t^{-2},t,p,w)u^3t^3
&= t^3w^3u^3 + u^2t^2w^2 p  + p^2 + p^2t^2u + 8ut^2p + 2ut^2 + 2u^2t^4 \\
& = \mathrm{PW}_{Y^\circ}(u,t,p,w).
\end{align*}

\section{P=W and SYZ}\label{sect:hkkp}

In this section, we will describe roughly how the mirror P=W conjecture can be seen as a result of SYZ mirror symmetry for log Calabi--Yau manifolds and a conjectural description of the behaviour of special Lagrangian torus fibrations on log Calabi--Yau manifolds. In this section, we will restrict ourselves to the case where $D$ is a smooth divisor in $X$. Remind also that we consider (co)homology with coefficients in $\mathbb{C}$.

We let $\omega$ be a symplectic form on $X$, and we let $\sigma$ be a section of the anticanonical bundle of $X$ so that $\sigma$ has simple vanishing along $D$. Then $\Omega = \sigma^{-1}$ is a nonvanishing section of the canonical bundle on $U$ which extends to a meromorphic section of the canonical bundle of~$X$ with a simple pole along $D$. In \cite{aur1}, Auroux studies special Lagrangian torus fibrations on $U= X \setminus D$. A torus $L$ embedded in $U$ by a map $\iota$ is called {\it special Lagrangian} if $\iota_*T_{L}$ is a maximal isotropic subspace of $T_U$ with respect to $\omega$ at each point in $L$, and $\mathrm{Im}(\Omega)|_L=0$. A fibration $\pi: U \rightarrow B$ is called a \emph{special Lagrangian torus fibration} if its smooth fibers are special Lagrangian tori. The SYZ mirror of $U$ is constructed by taking the moduli space $U^\vee$ of special Lagrangian torus fibers $L$
equipped with flat $U(1)$-bundles (which we indicate by the corresponding flat connections $\nabla$) up to gauge equivalence. Alternately, $U^\vee$ is the dual torus fibration over $B$, away from the critical locus.

In the case where $U$ is log Calabi--Yau, Auroux makes the following conjecture \cite[Conjecture 7.3]{aur1} about how special Lagrangian tori behave near the boundary of $B$.
\begin{conjecture}
The special Lagrangian torus fibers of $\pi$ near the divisor $D$ are $S^1$-bundles over a special Lagrangian torus fiber of a special Lagrangian torus fibration on $D$.
\end{conjecture}
More precisely, a tubular neighborhood $\mathscr{U}_D$ of $D$ in $X$ is symplectomorphic to a neighborhood of the zero section of the normal bundle of $D$ in $X$, hence $\mathscr{U}_D \setminus D$  is foliated by $S^1$-bundles over~$D$. This is an attractive conjecture for us because the punctured tubular neighborhood~$\mathscr{U}_D \setminus D$, or, more precisely, $\partial \overline{\mathscr{U}}_D$, plays a distinguished role in computing the weight filtration in the mixed Hodge structure on $H^i(U)$.

Let us recall that in the case where $D$ is smooth, $H^i(U)$ has only two nonzero weight graded pieces, $\Gr_W^iH^i(U)$ and $\Gr_W^{i+1}H^i(U)$, which are determined by the long exact sequence in cohomology
\[
\dots \longrightarrow H^{i-2}(D) \longrightarrow H^i(X) \longrightarrow H^i(U) \xrightarrow{\mathrm{res}} H^{i-1}(D) \longrightarrow H^{i+1}(X) \longrightarrow \dots.
\]
The map res is called the residue map. Particularly $\mathrm{Gr}^W_iH^i(U)\cong \mathrm{im}(H^i(X) \rightarrow H^{i}(U)) $, and~$\Gr^W_{i+1}H^i(U) =\ker(H^i(U) \rightarrow H^{i-1}(D))$. The weight filtration in homology is determined by the dual of this sequence, so $W_{-i-1}H_i(U)$ is the image of the homological dual to res. According to \cite[pp. 104]{pdelang}, the dual of the residue map is the ``tube over cycles map'', defined as follows. If $g\colon \partial \overline{\mathscr{U}}_D \rightarrow D$ is the natural map induced by the identification of $\mathscr{U}_D$ with a neighborhood of the identity in $N_{D/X}$, and $j\colon \partial \overline{\mathscr{U}}_D \rightarrow U$ is the natural embedding, then the tube over cycles map assigns to a cycle $c$ in $D$ the cycle $[j(g^{-1}(c))]$. In other words, $W_{-i-1}H_i(U)$ is spanned by cycles which are supported on $S^1$-fibrations over cycles in $D$.

Furthermore, Auroux gives us an idea of how such cycles behave under SYZ mirror symmetry.  We expect that the special Lagrangian torus fibers of $\pi$ are (close to) $S^1$-bundles over Lagrangian tori in $D$, and the $S^1$-fibers are leaves of the foliation mentioned above. In other words, as the special Lagrangian torus fiber $L$ gets closer and closer to $D$, the $S^1$-fiber should contract to a point. Therefore, there is a disc $\delta$ in $X$ bounded by $L$, and which intersects $D$ in a single point. One constructs a local holomorphic function on $U^\vee$ from $\delta$ defined as
\[
z_\delta(L,\nabla) = \exp\left( -\int_\delta \omega \right) \mathrm{hol}_\nabla(\partial \delta).
\]
This function is well-defined for pairs $(L,\nabla)$ where $L$ is near $D$, and should be a good approximation of the superpotential function $w$ on $U^\vee$ near the boundary of $B$.

Our goal now is to determine how a cycle in $W_{-i-1}H_i(U)$ behaves under torus duality. We will look at a very simple example which we believe exhibits many phenomena which appear in greater generality. One caveat here is that the homological point of view is likely not correct (or at the very least will be difficult to work with) in greater generality.

\begin{example}
Let $X$ be a rational elliptic surface with elliptic fibration $g\colon X \rightarrow \mathbb{P}^1$ and let $E$ be a smooth fiber of $g$.
There is expected to be a special Lagrangian torus fibration on $U = X \setminus E$, which we denote $\pi_U\colon U \rightarrow B$, whose Lagrangian fibers near $E$ look like $S^1$-bundles over special Lagrangian tori in $E$. In other words, we expect that near the boundary of the base $B$ there is a region diffeomorphic to $B_\infty = (0,1) \times S^1$, so that $\mathscr{U}_E$ is the preimage of $B_\infty$ under $\pi$ and the restriction of $\pi$ to $\mathscr{U}_E$ is the composition of an $S^1$-bundle over $(0,1) \times E$ and $\mathrm{id} \times \pi_E$, where~$\pi_E :E \rightarrow S^1$ is a special Lagrangian torus fibration on $E$. In other words, we have a morphism of torus bundles over $B_\infty$
\[
\mathscr{U}_E \longrightarrow (0,1) \times E.
\]
In fact, since we have chosen $E$ and $X$ so that $N_{E/X}$ is trivial, we have that $\mathscr{U}_E \cong (0,1) \times E \times S^1$.

Our task now is to identify homological cycles in $E$ spanning $H_1(E)$, determine their lifts in~$\partial \overline{\mathscr{U}}_E$, and see how they behave under torus duality. We may choose two copies of $S^1$ inside of $E$ which span $H_1(E)$, which are the fiber $F_E$ of $\pi_E$ and a section $S_E$ of $\pi_E$. The cycles span\-ning~$W_{-3}H_2(U)$ are then a fiber $F$ of $\pi_U$ restricted to $\partial \overline{\mathscr{U}}_E$ and an $S^1$-bundle over a section of $\pi_E$. Dualizing $\pi_U$ we dualize $\pi_E$,
hence determining the cycles in the dual torus fibration which are dual to $F_E$ and $S_E$ is straightforward.

The cycles which are dual to $F_E$ and $S_E$ are, respectively, a point in the dual torus fibration, and $p  \times E^\vee \times 1 \subset (0,1] \times E^\vee \times S^1 \cong\overline{\mathscr{U}}_E ^\vee$ for some point $p$ in $S^1$. Projection of $\overline{\mathscr{U}}_E^\vee \cong (0,1] \times E^\vee \times S^1$ onto $[0,1] \times S^1$ is essentially the map $z_\delta$ described above, therefore, the duals of the cycles in~$H_2(U)$ which span $W_{-3}H_2(U)$ should correspond to cycles supported in a fiber of $z_\delta$.

Near the boundary, the map $z_\delta$ is expected to behave, essentially, like the superpotential $w$ on the dual Landau--Ginzburg model. Therefore, the cycles corresponding to $F_E$ and $S_E$ are the images of the pushforward maps
\[
H_0(w^{-1}(t)) \longrightarrow H_0(U^\vee), \qquad H_2(w^{-1}(t)) \longrightarrow H_2(U^\vee)
\]
for $w^{-1}(t)$ a smooth fiber of $w$. In this case, the mirror of $U$ is expected to be, up to deformation,~$U$ itself, and the superpotential map is expected to be the elliptic fibration map $w\colon U \rightarrow \mathbb{C}$ induced by the given elliptic fibration $f\colon X \rightarrow \mathbb{P}^1$. Therefore, the perverse Leray filtration on cohomology is exactly the hyperplane filtration with respect to the stratification $p \in \mathbb{P}^1$. Thus
\[
P_{0}H^0(U^\vee) = \ker(H^0(U^\vee) \rightarrow H^0(w^{-1}(t))),\quad  P_{2}H^2(U^\vee) = \ker(H^2(U^\vee) \rightarrow H^2(w^{-1}(t))).
\]
We may define a dual filtration on homology. In which case, we have
\[
P_0H_0(U^\vee) = \mathrm{im}(H^0(w^{-1}(t)) \rightarrow H_0(U^\vee)), \quad P_{-2}H_2(U^\vee) = \mathrm{im}(H_2(w^{-1}(t)) \rightarrow H_2(U^\vee)).
\]
Therefore we have naturally identified $W_{-3}H_2(U)$ with $\bigoplus_i  P_{-i}H_i(U^\vee)$.
\end{example}

\section{Future directions}
 It is expected that P=W  conjecture is a broad  phenomenon.
 The above observations suggest that in addition to the
P=W conjecture for log Calabi--Yau varieties, a version of the~P=W
conjecture can be proved in the cases of  moduli spaces of bundles and
moduli spaces of representations of quivers.

Consider the case of moduli spaces of $SU(n)$-bundles $M_C$ over a curve $C$ and their
cotangent bundles $M_{Higgs}(C)$, which are moduli spaces of Higgs
bundles. %\footnote{What are duals for the bundles for which we make the P=W conjecture?}

\begin{question}
Does classical P=W conjecture for $M_{Higgs}(C)$
relate to  P=W conjecture for moduli space of $SU(n)$-bundles $M_C$?
\end{question}

Similarly given a quiver $Q$ we can consider a version of P=W conjecture for a moduli space~$M_{Q}$ of representations
of $Q$ and its cotangent bundle --- Nakajima variety
$M_{Naka}(Q)$.

\begin{question}
Does classical P=W conjecture for $M_{Q}$
relate to  P=W conjecture for Nakajima variety $M_{Naka}(Q)$?
\end{question}

Many  low-dimensional examples mentioned above are
obtained via gluing cluster varieties. This suggests
a possible way of approaching the relation between classical and mirror P=W
conjectures --- checking them on  cluster varieties and then
investigating the behaviour under gluing.

\end{document}